\documentclass[reqno]{amsart}

\usepackage{hyperref}
\hypersetup{
    colorlinks=true,
    linkcolor=blue,
    urlcolor=cyan,
    citecolor=blue
    }
\usepackage{bm}
\usepackage{amsmath}
\usepackage{amssymb}
\usepackage{enumitem}
\newcommand{\cone}{\mathcal{C}}
\newcommand{\ten}{\mathbf}
\newcommand{\gt}{\bm}

\newcommand{\scalar}[1]{\langle#1\rangle}

\newcommand{\norm}[1]{\lvert#1\rvert}
\newcommand{\str}{\varepsilon}
\newcommand{\dd}{\mathrm{d}}

\DeclareMathOperator{\sgn}{sgn}

\renewcommand{\[}{\begin{equation}}
\renewcommand{\]}{\end{equation}}

\usepackage[sort&compress, numbers]{natbib}

\newcommand{\revise}[1]{#1}

\newcommand{\R}{\mathbb{R}}

\newtheorem{theorem}{Theorem}
\newtheorem{lemma}{Lemma}
\newtheorem{proposition}{Proposition}

\theoremstyle{definition}
\newtheorem{definition}{Definition}

\theoremstyle{remark}

\begin{document}
 
\title{Isometric deformations of surfaces of translation}

\author{Hussein Nassar}
\address{Department of Mechanical and Aerospace Engineering, University of Missouri, Columbia, MO 65211, USA}
\thanks{This work is funded by NSF CMMI CAREER award no. 2045881. The author has no conflicts of interest. Data availability is not applicable as no data was generated.}
\email{nassarh@missouri.edu}
%
\maketitle


\begin{abstract}
A \emph{surface of translation} is a sum $(u,v)\mapsto\gt\alpha(u)+\gt\beta(v)$ of two space curves: a \emph{path} $\gt\alpha$ and a \emph{profile} $\gt\beta$. A fundamental problem of differential geometry and shell theory is to determine the ways in which surfaces deform isometrically, i.e., by bending without stretching. Here, we explore how surfaces of translation bend. Existence conditions and closed-form expressions for special bendings of the infinitesimal and finite kinds are provided. In particular, all surfaces of translation admit a purely torsional infinitesimal bending. Surfaces of translation whose path and profile belong to an elliptic cone or to two planes but never to their intersection further admit a torsion-free infinitesimal bending. Should the planes be orthogonal, the infinitesimal bending can be integrated into a torsion-free (finite) bending. Surfaces of translation also admit a torsion-free bending if the path or profile has exactly two tangency directions. Throughout, smooth and piecewise smooth surfaces, i.e., surfaces with straight or curved creases, are invariably dealt with and some extra care is given to situations where the bendings cause new creases to emerge. 
\end{abstract}

\section{Introduction}\label{sec1}
A surface of translation is the Minkowski sum of two space curves: to construct it, take every point of one curve and add it, as a vector, to every point of the other curve. Alternatively, it is a surface swept by translating one curve, the profile, along another, the path;\footnote{This is CAD terminology. Fancier, classical, terms are ``generatrix'' and ``directrix''.} in this construction, path and profile play equivalent roles. Should the profile and the path be discrete, i.e., piecewise linear, their sum becomes polyhedral with parallelogram faces. Polyhedral surfaces of translation are sometimes referred to as ``zonoids'' and have found their way into design and architecture \cite{Hanegraaf1980, Hart2021}. Other examples of surfaces of translation include cylinders and, more surprisingly, helicoids. \revise{Cloth, modeled as a net of inextensible fibers, embraces in static equilibrium certain surfaces of translation; see, e.g.,~\cite{Pipkin1983, Pipkin1993}}. In structural engineering, the most recognizable surfaces of translation are the folded corrugations used in roof panels and in sandwich panels \cite{Buannic2003, Biancolini2005, Lebee2010}.

Speaking of folds, several origami tessellations, including the notorious ``Miura ori'', are periodic polyhedral surfaces of translation. Their appeal resides in the fact that they can deploy from densely-packed states into spread-out states that shade large areas,\footnote{See, e.g., the recently deployed sunshield of the JWST: \url{https://webb.nasa.gov/content/observatory/sunshield.html}} simply by folding and unfolding along crease lines. Remarkably, throughout the deployment, these origami tessellations remain surfaces of translation \cite{Nassar2022}. Investigating how this property carries over to other, smooth and piecewise smooth, surfaces of translation is a main goal of the present work.

More generally, compliant shell mechanisms are shell structures ornate with folds and corrugations that can ``morph'' in order to adapt to changes in environment and loading conditions, for locomotion purposes or for deployment purposes \cite{norman2009phd, Norman2008, Norman2009, Schenk2011a, Seffen2012}. Accordingly, they admit multiple, or even a continuum of, low-energy configurations wildly different in terms of spans, curvatures and overall shapes. Ideally, the low-energy configurations of a shell correspond to surfaces that are mutually isometric. Indeed, shell deformations due to bending and stretching contribute energies of different orders of magnitude. Their a priori ratio is proportional to $h^2/R^2$, where $h$ is the thickness of the shell and $R$ is a characteristic radius of curvature \cite{Landau1986}. Thus, in the limit of infinitely thin shells, low-energy configurations are stretch-free, i.e., they have the same metric as the natural state. This is why understanding how surfaces deform isometrically can be beneficial for the design and modeling purposes of compliant shell mechanisms.

Beyond these applications, the above energy argument\footnote{Rayleigh \cite{Rayleigh1894} appears to be the first to formulate the argument. See also the brief historical note in Love's Treatise \cite{Love1906}.} shows that a thin shell will favor isometric deformations if available. In that case stretching energy vanishes and can be discarded; the resulting model is that of a ``flexural shell''. By contrast, if the shell geometry, in conjunction with the boundary conditions, does not allow for isometric deformations, then stretching energy will always dominate bending energy. Accordingly neglecting the latter leads to the model of a ``membrane shell''. Rigorous asymptotic considerations suggest that there is in fact a hierarchy of shell models depending on the magnitude of the applied loads; see discussions in \cite{ciarlet2006} and~\cite{Harutyunyan2017}.

Evidently, isometric deformations are also interesting intrinsically as a purely geometric topic. Hereafter, a few results from the literature are recounted mainly for the benefit of the newcomer and the curious reader. But first, it will prove convenient to introduce a terminology that permits to distinguish different notions of isometric deformations. An \emph{isometry} is a deformation that preserves the lengths of curves. An isometry is \emph{trivial} if it is a Euclidean motion, i.e., the composition of a translation and a rotation; otherwise, it is \emph{non-trivial}. A non-trivial isometry is further called a \emph{warping}. A \emph{(finite) bending} is a (continuous) one-parameter family of isometries that starts with the identity and includes a warping. An \emph{infinitesimal bending} is a velocity field that preserves the lengths of curves up to first order in time, and that is not the velocity field of a Euclidean motion. If a surface does not admit a warping, then it is \emph{globally rigid}. If it does not admit a bending, then it is \emph{rigid}; otherwise it is \emph{flexible}. Last, if it does not admit an infinitesimal bending, then it is \emph{infinitesimally rigid}; otherwise it is \emph{infinitesimally flexible}.

The earliest rigidity result of significance is Cauchy's theorem: convex polyhedra with the same net are congruent. Equivalently, there exists no warping of a convex polyhedron that preserves convexity. Clearly then, there exists no bending that preserves convexity either. Later, Dehn showed that a convex polyhedron\footnote{There is an interesting generalization of Dehn's theorem to non-convex polyhedra with vertices in convex position \cite{Izmestiev2010}.} cannot bend infinitesimally except when it has flat vertices.\footnote{If a vertex is flat, then a displacement normal at the vertex and vanishing everywhere else provides a, somewhat uninteresting, infinitesimal bending.} This leaves open two questions: ($i$) whether a convex polyhedron with flat vertices can bend ``inwards'' into a non-convex polyhedron; and, ($ii$) whether non-convex polyhedra are rigid. Connelly answered both questions: he demonstrated that a convex polyhedron with flat vertices is rigid and therefore that all convex polyhedra are rigid; he also constructed the first flexible non-convex polyhedron. For precise statements, see \cite{Ivanova-Karatopraklieva1995, Connelly1993}.

In this context of polyhedral surfaces, it was implicitly understood that isometries are trivial face-wise; that is, each face moves like a rigid body. There is a stronger version of Cauchy's theorem, due to Alexandrov, that makes no a priori assumptions regarding the deformation of the faces. An even stronger theorem is due to Pogorelov: isometric compact convex surfaces are congruent  \cite{Connelly1993, Ivanova-Karatopraklieva1994}. Here, the convex surfaces need not be polyhedral or regular beyond whatever regularity is implied by convexity. Evidently there could be warpings that do not preserve convexity. Even more, by a result due to Kuiper, convex, and in fact all, $C^1$-smooth surfaces admit a $C^1$-smooth bending. Note that this bending leads to non-convex surfaces that are, loosely speaking, corrugated on an infinitely fine scale; see, e.g., the illustration in \cite{Borrelli2013}. In particular, the image surfaces are not $C^2$-smooth and lack the notion of Gaussian curvature. Hereafter, smooth surfaces and isometries are understood to be at least of class $C^2$.

Several results on the rigidity of smooth compact convex surfaces preceded Pogorelov's theorem. A first proof of infinitesimal rigidity\footnote{Just as for Dehn's theorem, infinitesimal rigidity holds except when it trivially does not, i.e., for surfaces with planar neighborhoods.} was attempted by Jellett\footnote{In the Treatise \cite{Love1906}, Love references Jellett's theorem in confirmation to the conclusion that a complete spherical shell is infinitesimally rigid. Interestingly, in doing so, Love misquotes the theorem and appears to believe that ``closed surfaces'' are infinitesimally rigid, be them convex or not.} \cite{Jellett1849} but the first rigorous proof is attributed to Liebmann \cite{Spivak1999a}. Later, Cohn-Vossen proved that isometric smooth compact convex surfaces are congruent. As a direct consequence, smooth compact convex surfaces are globally rigid since the sign of Gaussian curvature encodes convexity and is an isometric-invariant. Cohn-Vossen is also credited with the first non-trivial example of an infinitesimally flexible smooth compact (non-convex) surface. No flexible smooth compact surfaces appear to be known however.

There are other interesting rigidity, and flexibility, results based on convexity, e.g., for closed non-compact convex surfaces and for compact convex surfaces minus a small neighborhood; see \cite{Spivak1999a}. Other results concern rigidity relative to a curve, i.e., with the image of a given curve prescribed, and local rigidity, i.e., in the vicinity of a point. In fact, many classical results belong to this local category whereas global results in the spirit of ``geometry in the large'' came later.\footnote{Perhaps the issue with classical results is not that they are local per se, but that they seldom informed on how different local behaviors can be sewed together; see, e.g., Spivak's discussion of the classical ``classification'' of developable surfaces \cite[Chapter~5, Section~4]{Spivak1999}. See also \cite{Demaine2011}.} In that regard, polyhedral surfaces are somewhat exceptional: Cauchy's and Dehn's theorems came first, and local characterizations of flexibility, or ``foldability'' as it is called in the context of origami, have only been sought relatively recently \cite{Streinu2004, Abel2016, Izmestiev2017}.

Results of a general character aside, significant efforts went into the derivation and classification of bendings for specific classes of surfaces: ruled, constant curvature, Weingarten, etc. The case of surfaces of revolution is particularly interesting, in applications granted, but also for it provides conceptually-important counter-examples to the plausible, but wrong, idea that flexibility comes from the boundary. On one hand, Cohn-Vossen's surface, referenced above, is an infinitesimally flexible compact surface of revolution. On the other hand, there are surfaces of revolution with boundaries that are infinitesimally rigid. See, e.g.,~\cite{Audoly2010} for illustrations and for a tractable treatment.

It is to this effort of ``derivation and classification'' that the present paper contributes. Throughout, closed-form expressions for certain infinitesimal bendings and bendings are provided. Expressions for bendings of surfaces of translation were first obtained by Bianchi. Interestingly, it was at this occasion that Bianchi baptized them ``superficie di traslazione'' \cite{Bianchi1878}. Bianchi's bendings are available in the case where the path and profile are smooth graphs with the same ``$y$-axis'' and perpendicular ``$x$-axes''. Here, similar bendings are investigated, namely
\begin{enumerate}[label=($\roman*$)]
    \item bendings and infinitesimal bendings where path and profile belong to two perpendicular planes (Proposition~\ref{thm:infbianchi} and Theorem~\ref{thm:bianchi});
    \item bendings and infinitesimal bendings where the profile has exactly two slopes (Theorems~\ref{thm:infbndParticular} and~\ref{thm:koko});
    \item infinitesimal bendings where tangents to path and profile belong to an elliptic cone that potentially degenerates into two, non-orthogonal, planes (Theorem~\ref{thm:infbnd}),
\end{enumerate}
all in the context of piecewise smoothness encompassing smooth surfaces as well as surfaces with straight and curved creases.

\section{Preliminaries}
First of all, a few definitions and examples are due.

\begin{definition}
Let $\Omega\subset\R^n$ be a non-empty bounded open set and let $\bar{\Omega}$ be its closure. A mapping $\ten x: \bar{\Omega}\to\R^m$ is \emph{piecewise smooth} if it is continuous and if there exists a finite number of disjoint connected open sets $\Omega_i$ such that $\cup_i\bar\Omega_i=\bar\Omega$ and $\ten x$ is smooth over $\bar\Omega_i$. The $\bar\Omega_i$ are then called the \emph{pieces} of~$\ten x$.
\end{definition}

Thus, hereafter, continuity is implicitly presumed in piecewise smoothness. The restriction to compact domains of definition might not be necessary but certainly simplifies technical statements especially where integration is needed. In the same fashion, the following definition of a surface is not the usual one but is suitable for the present purposes.

\begin{definition}
A \emph{surface} is a piecewise smooth mapping
\[
    \begin{split}
        \ten x: \R^2\supset \bar{\Omega}&\to\R^3\\
        (u,v)&\mapsto \ten x(u,v)
    \end{split}
\]
such that the partial derivatives $\ten x_1\equiv\ten x_u$ and $\ten x_2\equiv\ten x_v$ are everywhere linearly independent.
\end{definition}

A surface can therefore exhibit at the boundary of pieces discontinuities in the tangent plane that can be understood as crease lines. The concept of ``pieces'' will be useful to keep track of where creases are or of where they could appear as a result of a deformation.

\begin{definition}
Let $I$, $I_i$, $J$ and $J_j$ be non-empty bounded open intervals such that $\bar I = \bigcup_i\bar I_i$ and $\bar J = \bigcup_j\bar J_j$. A surface $\ten x$ defined over $\bar{\Omega}=\bar I\times \bar J$ with pieces $\bar\Omega_{i,j}=\bar I_i\times \bar J_j$ and such that $\ten x(u,v)=\gt\alpha(u)+\gt\beta(v)$ is a \emph{surface of translation}. The space curves $\gt\alpha$ and $\gt\beta$ are the \emph{path} and the \emph{profile} of the surface.
\end{definition}

For instance, $\ten x(u,v) = (\cos(u)+\cos(v),\sin(u)+\sin(v),u+v)$ is a surface of translation that is part of a helicoid. Note that the path and the profile describe the same helix. It does not escape the mindful reader that $\ten x_u$ and $\ten x_v$ are equal for $u = v~[2\pi]$ and the domain $\bar\Omega$ must avoid these lines. It is worth recalling here that the helicoid is a minimal surface. Other minimal surfaces of translation are Scherk's surfaces \cite{Gray2006}.

A surface of translation is solution to $\ten x_{uv}=\ten 0$. In particular, its second fundamental form is diagonal and the coordinate curves are said to be \emph{conjugate}. For $\gt\alpha$ and $\gt\beta$ piecewise linear, the surface of translation is polyhedral. The rectangular pieces of $\ten x$ then get mapped to parallelograms. In the language of discrete differential geometry, a quadrilateral surface whose quadrilaterals are planar is a \emph{discrete conjugate net} \cite{Schief2008}. Polyhedral surfaces of translation are therefore particular discrete conjugate nets. Examples include the ``Miura ori'' and ``eggbox'' patterns.

\section{Infinitesimal bendings}
A small deformation of a surface $\ten x$ can be described by a velocity field $\dot{\ten x}$. The deformed surface is then $\ten x + t\dot{\ten x}$ where time $t$ is understood to be a small scalar. The strains brought by the deformation, to leading order in $t$, are therefore
\[
    \str_{\mu\nu}=\frac{\scalar{\dot{\ten x}_\mu,\ten x_\nu}+\scalar{\dot{\ten x}_\nu,\ten x_\mu}}{2},\quad (\mu,\nu)\in\{1, 2\}^2.
\]
For these strains to vanish, $\dot{\ten x}$ must move the basis $(\ten x_u,\ten x_v)$ as a rigid body, i.e., there must exist a field of infinitesimal rotations $\ten w$ such that $\dot{\ten x}_\mu= \ten w\wedge\ten x_\mu$. Should $\ten w$ be constant, then the whole surface undergoes the same infinitesimal rotation and there is no actual ``bending''. This motivates the following definitions.

\begin{definition}
A \emph{Euclidean velocity} of $\ten x$ is a mapping $\dot{\ten x}=\ten v + \ten w\wedge\ten x$ where $\ten v$ and $\ten w$ are constants.
\end{definition}

\begin{definition}
An \emph{infinitesimal isometry} of a surface $\ten x$ is a piecewise smooth mapping $\dot{\ten x}$ with the same pieces as $\ten x$ such that
\[
\scalar{\dot{\ten x}_u,\ten x_u}=\scalar{\dot{\ten x}_v,\ten x_v} = \scalar{\dot{\ten x}_u,\ten x_v} + \scalar{\dot{\ten x}_v,\ten x_u}=0.
\]
If $\dot{\ten x}$ is not a Euclidean velocity then it is an \emph{infinitesimal bending}.
\end{definition}

\begin{theorem}\label{thm:uni}
All surfaces of translation are infinitesimally flexible. A ``universal'' infinitesimal bending is
\[
    \dot{\ten x}(u,v) = \gt\alpha(u)\wedge\gt\beta(v) + \int^u\gt\alpha\wedge\gt\alpha' -  \int^v\gt\beta \wedge\gt\beta'.
\]
\end{theorem}
\begin{proof}
The candidate $\dot{\ten x}$ is piecewise smooth since any jumps in $\gt\alpha'$ and $\gt\beta'$ are overcome by integration. The partial derivatives $\dot{\ten x}_u = (\gt\alpha-\gt\beta)\wedge\gt\alpha'$ and $\dot{\ten x}_v = (\gt\alpha-\gt\beta)\wedge\gt\beta'$ describe an infinitesimal rotation. Last, $\dot{\ten x}_{uv}=\gt\alpha'\wedge\gt\beta'\neq\ten 0$ implies that $\dot{\ten x}$ is not a Euclidean velocity.
\end{proof}

The provided infinitesimal bending twists the coordinate curves without ``bending'' them; in other words, it is a pure torsion. This is most easily seen for the plane $\ten x(u,v)=(u,v,0)$ for which $\dot{\ten x} = (0,0,uv)$. More generally, the velocity of the unit normal $\ten n$ is $\dot{\ten n} = (\gt\alpha-\gt\beta)\wedge\ten n$ since $\gt\alpha-\gt\beta$ is the infinitesimal rotation of the infinitesimal bending. Thus, letting
\[
    \dot e = \scalar{\dot{\ten x}_{uu},\ten n}+\scalar{\ten x_{uu},\dot{\ten n}},\quad
    \dot g = \scalar{\dot{\ten x}_{vv},\ten n}+\scalar{\ten x_{vv},\dot{\ten n}},\quad
    \dot f = \scalar{\dot{\ten x}_{uv},\ten n},
\]
be the changes inflicted to the coefficients of the second fundamental form, it is easily seen that
    \[
        \dot e = \dot g = 0, \quad
        \dot f =\norm{\gt\alpha'\wedge\gt\beta'}.
    \]
In fact, solving the linearized Gauss-Codazzi-Mainardi equations under the assumption $\dot e=\dot g = 0$ is what lead to the discovery of $\dot{\ten x}$ in the first place.

One might object to the ``universality'' of the above infinitesimal bending for, in the case of polyhedral surfaces, it does not preserve the planarity of the faces. But this should not be concerning in a setting where one of the goals is to invariably treat surfaces with, say, curved and straight creases. A more concerning objection would be as follows: suppose $\gt\beta$ describes a closed curve over an interval $[v_1,v_2]$, i.e., such that $\gt\beta(v_1)=\gt\beta(v_2)$ as it would for a ``tube of translation'' for instance. Then, it would be natural to ask of an infinitesimal bending to satisfy a similar constraint of \emph{closure}, namely, $\dot{\ten x}(u,v_1)=\dot{\ten x}(u,v_2)$. Clearly, the above ``universal'' infinitesimal bending does not satisfy said constraint and produces a ``dislocation'' with a ``Burger's vector''
\[
    \ten b = -\int_{v_1}^{v_2} \gt\beta\wedge\gt\beta'.
\]
The infinitesimal bendings described below are more hopeful in these regards but are not as universal.

\begin{proposition}\label{thm:infbianchi}
Let $\ten x = \gt\alpha + \gt\beta$ be a surface of translation where $\gt\alpha$ and $\gt\beta$ are not both straight lines. Suppose further that $\gt\alpha$ and $\gt\beta$ belong to two perpendicular planes such that $\gt\alpha'$ and $\gt\beta'$ are never in each other's plane. Then, up to a change of basis of $\R^3$,
\[
    \gt\alpha'=(x_\alpha,0,z_\alpha),\quad\gt\beta'=(0,y_\beta,z_\beta), \quad x_\alpha \neq 0, \quad y_\beta \neq 0,
\]
and
\[
    \dot{\ten x} = \int^u \dot{\ten x}_u + \int^v\dot{\ten x}_v
\]
where
\[
    \dot{\ten x}_u = \left(\frac{z^2_\alpha}{x_\alpha},0,-z_\alpha\right),\quad
    \dot{\ten x}_v = \left(0,-\frac{z^2_\beta}{y_\beta},z_\beta\right),\quad
\]
defines an infinitesimal bending.
\end{proposition}
\begin{proof}
    Let $P$ and $Q$ be the planes in which lie $\gt\alpha'$ and $\gt\beta'$, respectively. Let $\ten p\in Q^\perp\subset P$, $\ten q\in P^\perp\subset Q$ and $\ten r\in P\cap Q$ be unit vectors. Clearly, $(\ten p,\ten q,\ten r)$ is an orthonormal basis of $\R^3$ with $\scalar{\gt\alpha',\ten q}=\scalar{\gt\beta',\ten p} = 0$ as required. Furthermore, $x_\alpha\equiv \scalar{\gt\alpha',\ten p}\neq 0$ since otherwise $\gt\alpha'$ belongs to $Q$; similarly $y_\beta\neq 0$. Thus, $\dot{\ten x}$ is well-defined and is easily checked to be an infinitesimal isometry. Last, $\dot{\ten x}$ is not a Euclidean velocity. Indeed, suppose it was; then the infinitesimal rotation
    \[
        \ten w = \left(\frac{z_\beta}{y_\beta},\frac{z_\alpha}{x_\alpha}, \frac{z_\alpha z_\beta}{x_\alpha y_\beta}\right)
    \]
    is constant, and so are $z_\beta/y_\beta$ and $z_\alpha/x_\alpha$ meaning that both $\gt\alpha$ and $\gt\beta$ are straight lines.
\end{proof}

Perhaps the most pressing question is how does one come up with such an infinitesimal bending. Here is one way: it is known that the infinitesimal bendings of a smooth graph $(u,v,f(u,v))$ correspond to solutions $g$ of the linear PDE
\[\label{eq:infbendgraph}
    f_{vv}g_{uu} - 2f_{uv}g_{uv} + f_{uu}g_{vv} = 0.
\]
When the graph is a surface of translation, i.e., $f(u,v)=\alpha(u)+\beta(v)$, the equation simplifies into
\[
    \beta''g_{uu} + \alpha''g_{vv} = 0.
\]
The symmetry then immediately inspires the solution $g(u,v)=\alpha(u)-\beta(v)$. Working the way back to $\dot{\ten x}$ leads to the infinitesimal bending provided above.

The above infinitesimal bending is pleasant in many ways: it preserves conjugacy of the coordinate curves, i.e., it is free of torsion. Indeed, $\ten x +t\dot{\ten x}$ remains a surface of translation. Moreover, $\dot{\ten x}_u$, $\dot{\ten x}_v$ and $\ten w$ are constant wherever $\gt\alpha'$ and $\gt\beta'$ are constant. Hence, if $\ten x$ is polyhedral, then so is $\ten x+t\dot{\ten x}$. Last, if $\gt\alpha'$ and $\gt\beta'$ are periodic, then $\dot{\ten x}_u$ and $\dot{\ten x}_v$ are periodic with the same respective periods. There is a ``trick'', or two, that allows to construct similar infinitesimal bendings for surfaces of translation where $\gt\alpha'$ and $\gt\beta'$ belong to two, not necessarily orthogonal, planes or, more generally in a sense, to a single elliptic cone. On one hand, infinitesimal flexibility is a linear invariant: if $\dot{\ten x}$ is an infinitesimal bending of $\ten x$, then $\ten A^{-T}\dot{\ten x}$ is an infinitesimal bending of $\ten A\ten x$ where $\ten A$ is any invertible matrix. On the other hand, if $\dot{\ten x}$ is an infinitesimal bending of $\ten x$, then it is also an infinitesimal bending of $\ten x + \ten w\wedge\dot{\ten x}$ for any constant $\ten w$. Combining these two transformations leads to the following.

\begin{theorem}\label{thm:infbnd}
    Let $\ten x = \gt\alpha + \gt\beta$ be a surface of translation where $\gt\alpha$ and $\gt\beta$ are not both straight lines. Suppose further that $\gt\alpha'$ and $\gt\beta'$ belong to an elliptic cone $\cone$. If
    \begin{itemize}
        \item $\cone$ is non-degenerate; or
        \item $\cone$ is the union of two planes and $\gt\alpha'$ and $\gt\beta'$ each belong to one plane and never to the other plane,
    \end{itemize}
    then $\ten x$ admits an infinitesimal bending that preserves conjugacy, planarity and periodicity.
\end{theorem}
\begin{proof}
    Let $M$ be a plane of mirror symmetry such that $\cone\cap M$ contains two lines; let $\ten p$ and $\ten q$ be unit vectors along these lines and let $\ten r$ be a unit normal to $M$. Then, there exists $s\in\R$ such that
    \[
        \cone = \{x\ten p+y\ten q+z\ten r,\quad xy = sz^2\}.
    \]
    Let
    \[
        \gt\alpha' = x_\alpha\ten p + y_\alpha\ten q + z_\alpha\ten r,\quad
        \gt\beta' = x_\beta\ten p + y_\beta\ten q + z_\beta\ten r,
    \]
    and let $(\ten p^*,\ten q^*,\ten r^*)$ be the basis dual to $(\ten p,\ten q,\ten r)$. Then, by hypothesis,
    \begin{itemize}
        \item either $\cone$ is non-degenerate and $s\neq 0$, in which case let
    \[
        \dot{\ten x}_u = \frac{y_\alpha}{s}\ten p^* - z_\alpha\ten r^*,\quad
        \dot{\ten x}_v = -\frac{x_\beta}{s}\ten q^* + z_\beta\ten r^*;
    \]
    \item
    or $\cone$ is degenerate, $s=0$ but $x_\alpha\neq 0$ and $y_\beta\neq 0$ in which case let
    \[
        \dot{\ten x}_u = \frac{z_\alpha^2}{x_\alpha}\ten p^* - z_\alpha\ten r^*,\quad
        \dot{\ten x}_v = -\frac{z_\beta^2}{y_\beta}\ten q^* + z_\beta\ten r^*.
    \]
    \end{itemize}
    Either way, $\dot{\ten x}$ is an infinitesimal bending that preserves conjugacy, planarity and periodicity.
\end{proof}

The above infinitesimal bending corresponds to the infinitesimal rotation
\[
\ten w = \frac{1}{1 - s \frac{z_\alpha z_\beta}{x_\alpha y_\beta}}
\left(
            \frac{z_\beta}{y_\beta}\ten p
            +\frac{z_\alpha}{x_\alpha}\ten q
            + \frac{z_\alpha z_\beta}{x_\alpha y_\beta}\ten r
            \right).
\]
Note that $\ten w$ belongs either to a hyperbolic paraboloid ($s=0$) or to a one-sheeted hyperboloid ($s\neq 0$), both of which are doubly ruled surfaces. This correspondence between conjugacy-preserving infinitesimal bendings of surfaces of translation and doubly ruled surfaces can be foreseen as follows. Suppose $\ten x$ is smooth. By definition, $\dot{\ten x}_\mu = \ten w\wedge\ten x_\mu$. Then, integrability of $\dot{\ten x}_\mu$ implies
\[
    \dot{\ten x}_{uv} - \ten w\wedge\ten x_{u v}= \ten w_u\wedge\ten x_v = \ten w_v\wedge\ten x_u.
\]
Assuming, $\ten x$ and $\ten x+t \dot{\ten x}$ are both surfaces of translation, it comes that
\[
    \ten w_u\wedge\ten x_v = \ten w_v\wedge\ten x_u = 0.
\]
That is: there exist two functions $a$ and $b$ such that
\[
    \ten w_u(u,v) = a(u,v)\ten x_v(v), \quad
    \ten w_v(u,v) = b(u,v)\ten x_u(u),
\]
meaning that a point following the coordinate lines of $\ten w$, along $\ten w_u$ or $\ten w_v$, moves in a straight line. Thus, the surface to which $\ten w$ belongs must be doubly ruled. This property was observed by Smith \cite{Smith1905} as he attempted to characterize surfaces of translation with conjugacy-preserving bendings rather than infinitesimal bendings. More generally, Bianchi calls \emph{associate} two surfaces with parallel tangent planes such that a conjugate net on one corresponds to an asymptotic net on the other \cite{Eisenhart1909}. Here, $\ten x$ and $\ten w$ are associate. Indeed, on one hand, $\ten w_u\wedge\ten w_v$ is parallel to $\ten x_u\wedge\ten x_v$ so that the planes tangent to $\ten w$ and $\ten x$, for equal $(u,v)$, are parallel. On the other hand, both $\ten w_{uu}=a_u\ten x_v$ and $\ten w_{vv}=b_v\ten x_u$ are tangent vectors so that the coordinate curves on $\ten w$ form an asymptotic net whereas they formed a conjugate net on $\ten x$. The study of associate surfaces, their existence conditions in particular, becomes synonymous to the study of certain classes of infinitesimal bendings. For an $\ten x$ that is piecewise smooth, $\ten w$ is not even continuous and the correspondence becomes more involved. See, e.g., \cite{Schief2008} in the polyhedral setting.

Theorem~\ref{thm:infbnd} applies to an $\ten x$ that describes the lateral surface of a triangular prism and thus provides an infinitesimal bending $\dot{\ten x}$ that preserves planarity. Of course, in this case, $\dot{\ten x}$ does not preserve closure. In general, should $\cone$ be degenerate and $\gt\beta$ describe a closed plane curve over $[v_1,v_2]$, then $\dot{\ten x}$ produces a dislocation of Burger's vector
\[
    \ten b = -\int_{v_1}^{v_2} \frac{z_\beta^2}{y_\beta}\ten q^*.
\]
Remarkably, if $\cone$ is non-degenerate, then $x_\beta y_\beta=s z_\beta^2$ enforces $\ten b=\ten 0$ by the closure of $\gt\beta$. In that case, $\dot{\ten x}$ preserves closure as well. There are other cases where closure is preserved, e.g., whenever the closed curve is centrosymmetric.

There is a particular case of Theorem~\ref{thm:infbnd} that is worth highlighting as it admits an interesting generalization that goes beyond the theorem itself. Recall first that checking whether a finite set of lines belong to an elliptic cone amounts to checking whether a certain linear system of equations in 6 unknowns admit non-trivial solutions, the unknowns being the coefficients of a homogeneous polynomial of degree 2 in 3 indeterminates.\footnote{I.e, $P(x,y,z)=ax^2+by^2+cz^2+dyz+ezx+fxy$.} In particular, 5 lines always belong to an elliptic cone. Hence, Theorem~\ref{thm:infbnd} potentially applies to surfaces of translation $\ten x$ where, say, $\gt\alpha$ has at most three distinct slopes and $\gt\beta$ has at most 2 distinct slopes. Now here is the generalization.

\begin{theorem}\label{thm:infbndParticular}
Let $\ten x = \gt\alpha + \gt\beta$ be a surface of translation where $\gt\beta$ has exactly two slopes and $\gt\alpha'$ never belongs to the plane of $\gt\beta'$. Then, $\ten x$ admits an infinitesimal bending that preserves conjugacy, planarity and periodicity.
\end{theorem}
Note that since $\gt\beta$ has two slopes, its tangent $\gt\beta'$ lies in a plane, namely the span of the two slopes.
\begin{proof}
Let $\gt\beta_1$ and $\gt\beta_2$ be the two unit slopes of $\gt\beta$ so that 
\[
    \gt\beta' = \chi_1\scalar{\gt\beta',\gt\beta_1}\gt\beta_1 + \chi_2\scalar{\gt\beta',\gt\beta_2}\gt\beta_2
\]
where the $\chi_i$'s are indicator functions and let
\[
\dot{\gt\beta}_1=(\gt\beta_1\wedge\gt\beta_2)\wedge\gt\beta_1,
\quad\dot{\gt\beta}_2=-(\gt\beta_1\wedge\gt\beta_2)\wedge\gt\beta_2.
\]
Then, since $\scalar{\gt\alpha',\gt\beta_1\wedge\gt\beta_2}\neq 0$ by hypothesis,
\[
    \begin{split}
        \dot{\ten x}_v &= \chi_1\scalar{\gt\beta',\gt\beta_1}\dot{\gt\beta}_1 + \chi_2\scalar{\gt\beta',\gt\beta_2}\dot{\gt\beta}_2,\\
        \dot{\ten x}_u &= \left(\scalar{\gt\alpha',\dot{\gt\beta}_2}\gt\beta_1-\scalar{\gt\alpha',\dot{\gt\beta}_1}\gt\beta_2\right)
        \wedge\frac{\gt\alpha'}{\scalar{\gt\alpha',\gt\beta_1\wedge\gt\beta_2}},
    \end{split}
\]
define an infinitesimal bending $\dot{\ten x}$ that preserves conjugacy, planarity and periodicity.
\end{proof}

What Theorem~\ref{thm:infbndParticular} really is saying is that if $\gt\beta$ has two slopes only, then $\gt\alpha$ is capable of following its lead. This is why $\dot{\gt\beta}_1$ and $\dot{\gt\beta}_2$ are chosen at first, somewhat randomly, and $\dot{\ten x}_u$ is constructed afterwards so as to produce zero infinitesimal strains. In fact, $(\dot{\gt\beta}_1,\dot{\gt\beta}_2)$ belongs in principle to the 4D space $\gt\beta_1^\perp\times\gt\beta_2^\perp$, but 3 of its dimensions correspond to rigid body rotations of $\gt\beta$ that extend into rigid body rotations of $\ten x$; the remaining dimension extends into the $\dot{\ten x}$ provided above. Last, note that if $\gt\beta$ describes a closed curve, then $\dot{\ten x}$ preserves its closure. The same does not hold for $\gt\alpha$ in general; exceptions include cases where $\gt\alpha$ has at most two slopes as well.

Theorem~\ref{thm:uni} aside, cases where, say, $\gt\alpha$ is planar and $\gt\beta'$ is parallel to the plane of $\gt\alpha$ have been avoided. Should that occur, two possibilities present themselves. On one hand, suppose $\gt\beta'$ is parallel to the plane of $\gt\alpha$ at an isolated point $v$. The surface then contains a curve $\gt\alpha+\gt\beta(v)$ that completely lies in a tangent plane and at which, generically, the Gaussian curvature changes sign. Such a curve, with vanishing normal curvature, would be ``rigidifying'' since it is minimally bent and cannot bend further infinitesimally; see the discussion in~\cite{Audoly2010}. Note that theorem~\ref{thm:uni} remains indifferent: the infinitesimal bending it provides does not bend the curves of the conjugate net but twists them. On the other hand, suppose $\gt\beta'$ is parallel to the plane of $\gt\alpha$ over an open interval. The surface then contains an open portion of a plane and any normal displacement compactly supported within that portion provides a somewhat uninteresting infinitesimal bending. That said, the following is potentially of interest.

\begin{proposition}\label{thm:normal}
    Let $\ten x = \gt\alpha + \gt\beta$ be a surface of translation where $\gt\alpha$ is planar but not straight and where $\gt\beta'$ is parallel to the plane of $\gt\alpha$ over a union of pieces. Then, $\ten x$ admits an infinitesimal bending that preserves conjugacy, planarity and periodicity.
\end{proposition}
\begin{proof}
    Let $\ten N$ be a constant unit normal to the the plane of $\gt\alpha$. Let $J_j$ be the intervals where $\gt\beta'$ is parallel to the plane of $\gt\alpha$; let $\chi_j$ be their respective indicator functions and let $c_j$ be constants. Then,
    \[
        \dot{\ten x}_u=\ten 0,\quad
        \dot{\ten x}_v=\sum_j c_j\chi_j \ten N,
    \]
    provides the sought infinitesimal bending.
\end{proof}

It is implicitly understood here that the constants are chosen in a way that is compatible with periodicity where applicable. Evidently, the infinitesimal bending exploits the pieces boundaries to dissimulate or create creases. Interestingly, as soon as that happens, $\gt\beta'$ is ``kicked'' out of the plane of $\gt\alpha$ in a way that could allow for the other results of the present section to apply. More on that in the next section.

\section{Bendings}\label{sec:bendings}
In exceptional cases, it is possible to pursue the infinitesimal bendings of the previous section so as to produce a finite bending. These are cases where the infinitesimal bending happens to preserve whatever conditions gave it birth in the first place. But before stating the results, two definitions are due.

\begin{definition}
    Two surfaces $\ten x$ and $\ten y$ are \emph{isometric} if they have the same pieces and the same metric, namely,
    \[
        \scalar{\ten x_\mu,\ten x_\nu} = \scalar{\ten y_\mu, \ten y_\nu},\quad (\mu,\nu)\in\{1, 2\}^2.
    \]
    They are \emph{congruent} if there exist an orthogonal matrix $\ten R$ and a vector $\ten T$ such that $\ten y = \ten R\ten x +\ten T$. Otherwise, they are \emph{warpings} of one another.
\end{definition}

\begin{definition}
    Let $I$ be an interval and $\ten x$ be a surface. A continuous one-parameter family of surfaces $I\ni t\mapsto\ten y(t)$ is a \emph{bending} of $\ten x$ if $\ten y(t)$ and $\ten x$ are identical for some~$t$, warpings of one another for some $t$, and isometric for all $t$.
\end{definition}

Revisiting Proposition~\ref{thm:infbianchi}, provided $\gt\alpha$ and $\gt\beta$ are in perpendicular planes, it is seen that an infinitesimal bending exists and keeps $\gt\alpha$ and $\gt\beta$ in perpendicular planes. By integration, the following ``groomed'' version of Bianchi's original result~\cite{Bianchi1878} is obtained; see also~\cite{Izmestiev2023} for a recent treatment in the context of ``T-surfaces''.

\begin{theorem}\label{thm:bianchi}
    Let $\ten x=\gt\alpha + \gt\beta$ be a surface of translation where $\gt\alpha$ and $\gt\beta$ belong to two perpendicular planes and are not both straight lines. Suppose further that the domain where $\gt\alpha'$ and $\gt\beta'$ are in the intersection of their planes is either empty or a union of pieces and piece boundaries. Then, $\ten x$ admits a bending that preserves conjugacy, planarity and periodicity.
\end{theorem}
\begin{proof}
    For simplicity, reparametrize $\gt\alpha$ and $\gt\beta$ by arc length and let
    \[
        \gt\alpha'=(x_\alpha,0,z_\alpha),\quad
        \gt\beta'=(0,y_\beta,z_\beta).
    \]
    Note then that $\norm{z_\alpha}$ and $\norm{z_\beta}$ are both smaller than 1 and cannot both reach~1 since otherwise $\gt\alpha'(u)=\pm\gt\beta'(v)$ holds for some $(u,v)$. Thus, $I\equiv \mathopen]\max \norm{z_\alpha},1/\max\norm{z_\beta}\mathclose[$ is non-empty and $1\in\bar I$. Now let $\sgn(x_\alpha)$ be a piecewise constant function, with the same pieces as $\ten x$, that returns the sign of $x_\alpha$ where $x_\alpha\neq 0$ and $\pm 1$ otherwise. Define $\sgn(y_\beta)$ in the same fashion. Then,
    \[
        \begin{split}    
            \ten y_u(t) &= \left(\sgn(x_\alpha)\sqrt{1-z^2_\alpha/t^2},0,z_\alpha/t\right),\\
            \ten y_v(t) &= \left(0,\sgn(y_\beta)\sqrt{1-t^2z^2_\beta},t z_\beta\right),
        \end{split}
    \]
    defined over $\bar I$ is the sought bending.
\end{proof}

Here too, regularity is recovered by integration. The only difficulty resides in the fact that the points where $x_\alpha$ (or $y_\beta$) reaches $0$ could be regular at $t=1$ and become singular at $t\neq 1$. This means that $\ten y$, for $t\neq 1$, could exhibit some bending-induced creases; the difficulty is resolved by anticipating the emergence of these creases and by making suitable assumptions regarding the pieces of $\ten x$. On a similar note, it could be surprising that the above theorem is slightly more general than, say, Proposition~\ref{thm:infbianchi}, for it allows for the presence of the ``rigidifying'' curves discussed at the end of the previous section. However, the ``initial velocity'', i.e., $\dd\ten y/\dd t$ at $t=1$, is singular in such cases and does not define an infinitesimal bending.

It is worth noting that there is a degree of arbitrariness in the way $\sgn(x_\alpha)$ and $\sgn(y_\beta)$ were defined. Should different choices be possible, different bendings would be available, each distinguished by the set of creases it produces. Furthermore, different ``chunks'' of these different bendings could be stitched together at $t=1$ to provide other bendings. In cases where periodicity is to be preserved, only periodic choices of $\sgn(x_\alpha)$ and $\sgn(y_\beta)$ would be acceptable. Last, note that $\ten y(t)$ does not preserve closure in general; exceptions include cases where the closed curve is centrosymmetric.

Back to the thread of the present section, the other category of surfaces where the infinitesimal bending preserves its existence condition is the one in Theorem~\ref{thm:infbndParticular}. Before stating the second main theorem of the section, a few preliminaries are needed.
\begin{definition}
    Let $\ten x_1$, $\ten x_2$ and $\ten x_3$ be three unit vectors. The triplet $(\ten x_1,\ten x_2,\ten x_3)$ is said to be \emph{flat-folded} (resp. \emph{flat-unfolded}) if $\ten x_3 = a\ten x_1 + b\ten x_2$ with $ab\leq 0$ (resp. $ab\geq 0$) or if $\ten x_1=\ten x_2$ (resp. if $\ten x_1=-\ten x_2$).
\end{definition}

\begin{lemma}\label{lem:trineq}
    Let $\ten x_1$, $\ten x_2$ and $\ten x_3$ be three unit vectors and let $c_1=\scalar{\ten x_2,\ten x_3}$ and $s_1=\sqrt{1-c^2_1}$ and so on. Then,
    \[
        c_1c_2 - s_1s_2 \leq c_3 \leq c_1c_2 + s_1s_2,
    \]
    with equality to the left (resp. right) if and only if $(\ten x_1, \ten x_2, \ten x_3)$ is flat-unfolded (resp. flat-folded).
\end{lemma}
\begin{proof}
    Left as an exercise.
\end{proof}

\begin{theorem}\label{thm:koko}
    Let $\ten x = \gt\alpha + \gt\beta$ be a surface of translation where $\gt\beta$ has exactly two slopes $\gt\beta_1$ and $\gt\beta_2$ and $(\gt\beta_1,\gt\beta_2,\gt\alpha')$ is never flat-unfolded or never flat-folded. Suppose further that the domain where $(\gt\beta_1,\gt\beta_2,\gt\alpha')$ is flat-folded or flat-unfolded is either empty or a union of pieces and piece boundaries. Then, $\ten x$ admits a bending that preserves conjugacy, planarity and periodicity.
\end{theorem}
\begin{proof}
    Re-parametrize, for simplicity, $\gt\alpha$ and $\gt\beta$ by arc length. Let $\gt\beta_1$ and $\gt\beta_2$ be the two unit slopes of $\gt\beta$ so that 
    \[
        \gt\beta' = \chi_1\scalar{\gt\beta',\gt\beta_1}\gt\beta_1 + \chi_2\scalar{\gt\beta',\gt\beta_2}\gt\beta_2
    \]
    where the $\chi_i$'s are indicator functions. Let
    \[
        c_\mu = \scalar{\gt\alpha',\gt\beta_\mu},\quad
        s_\mu = \sqrt{1-c_\mu^2},\quad \mu = 1,2, \quad
        c = \scalar{\gt\beta_1,\gt\beta_2}.
    \]
    Then, by hypothesis and by Lemma~\ref{lem:trineq},
    \[
        \begin{split}
        \text{either}\quad\max(c_1c_2-s_1s_2)&\leq c< \min(c_1c_2+s_1s_2),\\
        \text{or}\quad\max(c_1c_2-s_1s_2)&< c \leq \min(c_1c_2+s_1s_2).
        \end{split}
    \]
    In any case, $I\equiv\mathopen]\max(c_1c_2-s_1s_2),\min(c_1c_2+s_1s_2)\mathclose[$ is non-empty and $c\in\bar I$. For $t\in\bar I$, let $\gt\beta_1(t)$ and $\gt\beta_2(t)$ be two unit vectors such that
    \[
        \gt\beta_1(c) = \gt\beta_1,\quad
        \gt\beta_2(c) = \gt\beta_2,\quad
        \scalar{\gt\beta_1(t),\gt\beta_2(t)} = t.
    \]
    The detail of this construction is not important. Let
    \[
        \begin{split}
        \ten a(t) &= c_1\gt\beta_2(t) - 
        c_2\gt\beta_1(t),\\
        \ten b(t) &= \gt\beta_1(t)\wedge\gt\beta_2(t),\\
        J &= \scalar{\gt \beta_1\wedge\gt \beta_2,\gt\alpha'},
        \end{split}
    \]
    and let $\sgn(J)$ be a piecewise constant function, with the same pieces as $\ten x$, that returns the sign of $J$ if $J\neq 0$ and $\pm 1$ otherwise. Then,
    \[
    \begin{split}
        J(t) &\equiv \sgn(J)\sqrt{\norm{\ten b(t)}^2-\norm{\ten a(t)}^2} \\
             & = \sgn(J)\sqrt{(c_1c_2+s_1s_2-t)(t-c_1c_2+s_1s_2)}
    \end{split}
    \]
    is well-defined for $t\in\bar I$. Finally,
     \[
        \begin{split}
        \ten y_v(t) &= \chi_1\scalar{\gt\beta',\gt\beta_1}\gt\beta_1(t) + \chi_2\scalar{\gt\beta',\gt\beta_2}\gt\beta_2(t),\\
        \ten y_u(t) & = \frac{\ten a(t)\wedge\ten b(t) + J(t)\ten b(t)}{\norm{\ten b(t)}^2}
        \end{split}
    \]
    provide the sought bending.
\end{proof}

The above theorem affords the same discussion as Theorem~\ref{thm:bianchi} regarding the emergence of new creases, its generality compared to Theorem~\ref{thm:infbndParticular}, and the arbitrariness in the definition of $\sgn(J)$. Note that $\ten y(t)$ preserves the closure of $\gt\beta$; it does not necessarily preserve the closure of $\gt\alpha$.

For discrete $\gt\alpha$, the existence part in Theorem~\ref{thm:koko} could have been obtained differently. Indeed, it is known \cite{Stachel2010} that a mesh of $3\times 3$ planar quadrilaterals\footnote{Also known as a quad-based Kokotsakis mesh.} of the ``translational'' type is flexible. Therein, being of the ``translational'' type is equivalent to being a polyhedral surface of translation where the path or the profile has exactly two slopes. Combine that with a theorem from~\cite{Schief2008} that states that a ``non-degenerate'' discrete conjugate net is flexible if and only if its $3\times 3$ complexes are flexible. The obtained result guarantees the existence of a bending for ``non-degenerate'' discrete conjugate nets of translation where the path or the profile has exactly two slopes. Theorem~\ref{thm:koko} generalizes this result to smooth and piecewise smooth settings and provides an expression for the bending.

\section{Conclusion}
To summarize, all surfaces of translation admit a purely torsional infinitesimal bending. Surfaces of translation whose path and profile belong to an elliptic cone or to two planes but never to their intersection further admit a torsion-free, i.e., conjugacy-preserving, infinitesimal bending. Should the planes be orthogonal, the infinitesimal bending can be integrated into a torsion-free (finite) bending. Surfaces of translation also admit a torsion-free bending if the path or the profile has exactly two tangency directions. These existence conditions, as well as the corresponding closed-form expressions, are equally valid for smooth and piecewise smooth surfaces, i.e., surfaces with straight or curved creases. It is not known if these conditions exhaust all possibilities or guarantee uniqueness.

\end{document}